\documentclass[11pt]{amsart}
\usepackage{amsrefs, hyperref}
\usepackage{nopageno}

\title{The user's guide project: giving experiential context to research papers}

\author{Cary Malkiewich}
\address{Department of Mathematics, Binghampton University, Binghamton, New York, USA}
\email{malkiewich@math.binghamton.edu}

\author{Mona Merling}
\address{Department of Mathematics, Johns Hopkins University, Maryland, USA}
\email{mmerling@math.jhu.edu}

\author{David White}
\address{Department of Mathematics and Computer Science \\ Denison University Granville, OH, USA}
\email{david.white@denison.edu}

\author{Franc Lucas Wolcott}
\address{Department of Mathematics, Lawrence University, Wisconsin, USA}
\email{luke.wolcott@gmail.com}

\author{Carolyn Yarnall}
\address{Department of Mathematics, University of Kentucky, USA}
\email{carolyn.yarnall@gmail.com}

\begin{document}

\begin{abstract}
This paper was written in 2015, and published in the {\em Journal of Humanistic Mathematics}. This paper announces the first issue (2015) of {\em Enchiridion: Mathematics User's Guides}, a project meant to produce peer-reviewed User's Guides as companions to published papers. These User's Guides are meant to explain the key insights and organizing principles in their companion papers, the metaphors and imagery used by the authors, the story of the development of the companion papers, and a colloquial summary appropriate for a non-mathematical audience. Examples of User's Guides can be found at \href{https://mathusersguides.com/}{https://mathusersguides.com/}.
\end{abstract}

\maketitle
\thispagestyle{empty}

\begin{quotation}
Nowhere in the sciences does one find as wide a gap as that between the written version of a mathematical result and the discourse that is required in order to understand the same result. -- \textit{Gian-Carlo Rota}~\cite{rota}\\
\end{quotation}

Have you ever read a math paper and wondered, how did they come up with this?  Why? What's really going on here?  What's the best way to think about this?

We are writing to announce the first issue of the User's Guide Project.  We, five topologists, have come together, each writing a user's guide to one of our published or soon-to-be-published research papers \cite{malkiewich,merling,white,wolcott,yarnall}.  We collaboratively group-peer-reviewed the guides, and they are available at \url{http://www.mathusersguides.com}.

A user's guide -- at the same time humanistic and technical -- is written to accompany a research article, providing further exposition and context for the results.  Our user's guides are composed of four topics.\\

\noindent \textit{Topic 1: Key insights and organizing principles.} 

What is the conceptual essence of the paper? What's really going on? What is the best -- most natural, most concise, most elegant -- overall conceptual framework for all the ideas in play?

This description differs from your average paper introduction in the following way.  In an introduction, we need to give a quick overview of the content and situate our results in a mathematical context, and then clearly announce our new results and contributions.  There is pressure to convince the reader why they should care and why they should read the paper.  In the user's guide there is no pressure to advertise; the ideas and results themselves are given central importance.\\

\noindent \textit{Topic 2: Metaphors and imagery.}

While developing the ideas  that go into a paper, a mathematician often constructs useful metaphors and mental imagery to aid in understanding.  These insightful tools for accessing the results, while they may be divulged in an off-hand comment during a lecture or when talking with a colleague~\cite{ThurstonMO}, are almost always omitted from published mathematical literature.

Each of the authors in the User's Guide Project were asked to articulate these subjective modes of reasoning.  What mental images do you have when you think about these ideas?  What conceptual metaphors do you use?  What's ``the right way to think about it"?\\

\noindent \textit{Topic 3: Story of the development.}

Where and when did these ideas and results arise?  What is the logistic and psychological context of this paper?  What was the process of arriving at these particular statements and proofs?

This section tells a story --  about the background work, the dead ends, the moments of insight, the relevant conversations with others.  Traditional papers aim to give the shortest logic path between A and B, but here we tell a little about  the actual meandering path taken to arrive at the results.\\

\noindent \textit{Topic 4: Colloquial summary.}

The target audience for our user's guides are students and mathematicians that have the specialized training, and the desire,  to read the original research papers.  We envision the source paper and the user's guide being read side-by-side.  But each user's guide also includes a short colloquial summary, aimed at a much broader, non-mathematical audience, in which we try to explain at least a small part of the results, using everyday metaphors.  

Here we also share broadly the core of what drove us to work in this area, including to our colleagues, family, and undergraduate students who are just starting out. What drew us to conduct this work and what other endeavors of human expression can it be related to? \\ \\

Our goal with the User's Guide Project is to make research mathematics more accessible.  We would like to help other topologists who wish to read and understand our papers, but more broadly we would like to break out of current research literature conventions, exploring more flexible, and hopefully more insightful, expositional styles.

As in~\cite{LW-JHM}, we ``would like to think of this as meta-data for the ideas and results, which attaches a bit of humanity to the objective representations and reasoning. Including this information closes the gap between the practice of mathematics and the artifacts of that practice... But it must be clear the intention is not to remove rigor, just augment it, and perhaps re-prioritize it."

We believe these user's guides will be of interest to the readers of the \textit{Journal of Humanistic Mathematics}.  The first five user's guides, collaboratively peer-reviewed by the five of us, will be available in completed form in Fall 2015 at the project website \url{http://www.mathusersguides.com}.  Before that time, you will find a partial collection of topics.

We are looking for future writers as well.  If you have written a conventional theorems-and-proofs paper and would like to write a user's guide for it, please contact Luke Wolcott at \texttt{luke.wolcott@gmail.com} or David White at \texttt{david.white@denison.edu}.  Individual topics are written and posted throughout the year, and then compiled in the summer and fall.


\end{document}